\documentclass[a4paper,12pt]{article}
\usepackage{setspace}
\usepackage{fullpage}
\usepackage{amsthm}
\usepackage{amssymb}
\usepackage{amsmath}
\usepackage{graphicx}
\usepackage{caption}
\usepackage{subcaption}
\usepackage{color}
\usepackage{float}
\theoremstyle{definition}
\newtheorem{definition}{Definition}[section]
\theoremstyle{remark}

\theoremstyle{assumption}

\theoremstyle{example}
\newtheorem{example}{Example}[section]
\theoremstyle{property}

\theoremstyle{theorem}

\theoremstyle{lemma}
\newtheorem{lemma}{Lemma}[section]
\theoremstyle{corollary}

\theoremstyle{result}
\newtheorem{result}{Result}[section]

\begin{document}

\title{New fractional differential inequalities with their implications to the stability analysis of fractional order systems}

\author{Bichitra Kumar Lenka\\Department of Mathematics and Statistics,\\ Indian Institute of Science Education and Research Kolkata,\\
Nadia-741246, West Bengal, India\\
$Email$: bichitra.lenka99@gmail.com}

\maketitle

\begin{abstract}
It is well known that the Leibniz rule for the integer derivative of order one does not hold for the fractional derivative case when the fractional order lies between 0 and 1. Thus it poses a great difficulty in the calculation of fractional derivative of given functions as well as in the analysis of fractional order systems. In this work, we develop a few fractional differential inequalities which involve the Caputo fractional derivative of the product of continuously differentiable functions. We establish some of their properties and propose a few propositions. We show that these inequalities play a very essential role in the Lyapunov  stability analysis of nonautonomous fractional order systems.
\end{abstract}

\section{Introduction}\label{Sec1}
Fractional calculus plays a great role in many areas of science, engineering and interdisciplinary subjects, and has a broad range of applications \cite{Podlubny1999,Kilbas2006,Petras2011,Kiryakova2014}.  Fractional operators (fractional derivative and/or integral) are global operators \cite{Podlubny1999,Kilbas2006}. Unlike the usual local derivative operators (integer order derivative operators), which have easy properties (e.g. Leibniz rule, Chain rule and so on), the properties of fractional derivative operators are not simple  (e.g. Leibniz rule, Chain rule and so on).  Indeed, fractional operators possess very rich and complex properties \cite{Podlubny1999,Kilbas2006}.  In \cite{Tarasov2013,VETarasov2016,Tarasov2016}, it has been shown that the Leibniz rule, Chain rule, etc., do not hold for fractional derivatives whenever the fractional order lies between $0$ and $1$. In fact, these rules do not hold for the Caputo fractional derivatives too. Thus, due to the global nature and the complicated properties of Caputo fractional derivative operator, in practice, it is a very difficult task for the estimation of fractional derivatives of the given functions.
 
Fractional differential inequalities (which involve fractional derivatives) exist in the theory of fractional order systems. Indeed, they play a very essential role in the qualitative theory of fractional order systems. It may be noted that stability is an important qualitative property of the solutions to fractional order systems. In the literature, there has been a growing interest in the investigation of the stability or asymptotic stability of the solutions to the autonomous and nonautonomous fractional order systems \cite{YLi2009,YLi2010,Delavari2012,Yu2013,Mermoud2015,BKL2019,Fernandez2017}. 
 
It may be noted that Lyapunov direct method \cite{Khalil2002}, fractional Lyapunov direct method \cite{YLi2009,YLi2010,Delavari2012,Yu2013,Mermoud2015}, distributed order Lyapunov direct method \cite{Fernandez2017} and fractional comparison method \cite{BKL2019} are powerful methods for the stability analysis of the fractional order systems. However, the application of these methods require some suitable candidate functions or time-varying  Lyapunov functions (which involve both state variables (dependent variables) and time (independent variable)) that are continuously differentiable, and the calculation of its fractional derivatives. Indeed, the application of these methods require fractional differential inequalities.  

Recently, in \cite{Camacho2014,Mermoud2015,Kirane2015,Kirane2017,Dai2017,
Fernandez2017b,Fernandez2018,BKLenka2019}, the authors have developed some interesting fractional differential inequalities which involve the Caputo fractional derivatives. In this paper, we present some new fractional differential inequalities which involve the Caputo fractional derivatives of the product of continuously differentiable functions. Then, we propose some equivalence results and also establish some of their interesting properties. By presenting a few interesting examples, we demonstrate the usefulness of some appropriate fractional differential inequalities and their importance in the application of  fractional Lyapunov direct method.

\section{Notations and definitions}
Let us denote by $\mathbb{N}$ be the set of natural numbers, $\mathbb{Z}^{+}$ the set of positive integers, $\mathbb{R}^{+}$ the set of positive real numbers, $\mathbb{R}$ the set of real numbers, $\mathbb{C}$ the set of complex numbers, $\mathfrak{R}(z)$ the real part of complex number $z$, $\mathbb{R}^{n}$ the $n$-dimensional Euclidean space.

\begin{definition}\cite{Podlubny1999}
The Gamma function is defined as
\begin{equation}
\Gamma(z)=\int_{0}^{\infty}t^{z-1}e^{-t}dt
\end{equation}
where $z\in \mathbb{C}$ and $\mathfrak{R}(z)>0$.
\end{definition}

\begin{definition}\cite{Podlubny1999,Kilbas2006}
The Riemann-Liouville fractional integral of order $\alpha$ of a continuous function $x:[t_0, T)\to \mathbb{R}$, $-\infty<t_0< T \leq \infty$ is defined as 
\begin{equation}
^{RL}\!D_{t_0,t}^{-\alpha}x(t)=\frac{1}{\Gamma(\alpha)}\int_{t_0}^{t}\left(t-\tau\right)^{\alpha-1}x(\tau)d\tau
\end{equation}
where $\alpha\in \mathbb{R}^{+}$.
\end{definition}

\begin{definition}\cite{Podlubny1999,Kilbas2006}
The Caputo fractional derivative of order $\alpha$ of a $n^{th}$ continuously differentiable function $x:[t_0, T)\to \mathbb{R}$, $-\infty<t_0< T \leq \infty$ is defined as 
\begin{equation}
^{C}\!D_{t_0,t}^{\alpha}x(t)=\left\lbrace{\begin{array}{cc}
{^{RL}\!D_{t_0,t}^{-(n-\alpha)} \left(\frac{d^{n}x(t)}{dt^{n}}\right)}, & \text{if}\;\; \alpha \in \left(n-1, n\right)  \\
\frac{d^{n}x(t)}{dt^{n}}, & \text{if}\;\; \alpha=n
\end{array}}\right.
\end{equation}
where $\alpha\in \mathbb{R}^{+}$ and $n\in \mathbb{Z}^{+}$.
\end{definition}

\section{Fractional differential inequalities}\label{Sec2}

\subsection{Main inequalities}

\begin{result}\emph{\cite{BKLenka2019}}\label{NR1}
Let $\phi:[t_0,\infty)\to \mathbb{R}$ be a monotonically decreasing and continuously differentiable function. Suppose $x:[t_0,\infty) \to \mathbb{R}$ is a non-negative and continuously differentiable function. Then, the inequality 
\begin{align}\label{Res1_eq}
^{C}\!D^{\alpha}_{t_0,t}\left\lbrace{\phi(t)x(t)}\right\rbrace \leq {\phi(t)} {^{C}\!D^{\alpha}_{t_0,t}x(t)}, \;\; \forall t\geq t_0, \;\; \forall \alpha \in (0, 1],
\end{align}
holds.
\end{result}

\begin{result}\emph{\cite{BKLenka2019}}\label{Nr1}
Let $\phi:[t_0,\infty)\to \mathbb{R}$ be a monotonically increasing and continuously differentiable function. Suppose $x:[t_0,\infty) \to \mathbb{R}$ is a non-negative and continuously differentiable function. Then, the inequality 
\begin{align}
^{C}\!D^{\alpha}_{t_0,t}\left\lbrace{\phi(t)x(t)}\right\rbrace \geq {\phi(t)} {^{C}\!D^{\alpha}_{t_0,t}x(t)}, \;\; \forall t\geq t_0, \;\; \forall \alpha \in (0, 1],
\end{align}
holds.
\end{result}

\begin{result}\label{Nr2}
Let $\phi:[t_0,\infty)\to \mathbb{R}$ be a monotonically decreasing and continuously differentiable function. Suppose $x:[t_0,\infty) \to \mathbb{R}$ is a non-negative  continuously differentiable function. Then,  $ \forall t\geq t_0$, $\forall \alpha \in (0, 1]$, the following inequality holds
\begin{align}\label{Nr2_eq1}
^{C}\!D^{\alpha}_{t_0,t}\left\lbrace{\phi^{2n+1}(t)x^{\beta}(t)}\right\rbrace \leq {\phi^{2n+1}(t)} {^{C}\!D^{\alpha}_{t_0,t}x^{\beta}(t)}, 
\end{align}
where  $n\in \mathbb{N}$ and the constant $\beta$ is a non-negative real number.
\end{result}

\begin{result}\label{NR3}
Let $\phi:[t_0,\infty)\to \mathbb{R}$ be  a non-negative, monotonically decreasing and continuously differentiable function. Suppose $x:[t_0,\infty) \to \mathbb{R}$ is a non-negative  continuously differentiable function. Then,  $ \forall t\geq t_0$, $\forall \alpha \in (0, 1]$, the following inequality holds
\begin{align}
^{C}\!D^{\alpha}_{t_0,t}\left\lbrace{\phi(t)x^{\beta}(t)}\right\rbrace \leq {\phi(t)} {^{C}\!D^{\alpha}_{t_0,t}x^{\beta}(t)}, 
\end{align}
where  the constant $\beta$ is a non-negative real number.
\end{result}
\begin{proof}
The proof follows from the Result~\ref{NR1}.
\end{proof}

\begin{result}\label{Nr3}
Let $\phi:[t_0,\infty)\to \mathbb{R}$ be  a non-negative, monotonically decreasing and continuously differentiable function. Suppose $x:[t_0,\infty) \to \mathbb{R}$ is a continuously differentiable function. Then,  $ \forall t\geq t_0$, $\forall \alpha \in (0, 1]$, the following inequality holds
\begin{align}
^{C}\!D^{\alpha}_{t_0,t}\left\lbrace{\phi(t)x^{\beta}(t)}\right\rbrace \leq {\phi(t)} {^{C}\!D^{\alpha}_{t_0,t}x^{\beta}(t)}, 
\end{align}
where  the constant $\beta=2k$, and $k\in \mathbb{N}$.
\end{result}
\begin{proof}
The proof directly follows from the Result~\ref{NR3}.
\end{proof}
\subsection{Main propositions}
\begin{result}\label{NR4}
Let $\phi:[t_0,\infty)\to \mathbb{R}$ be  a non-negative, monotonically decreasing and continuously differentiable function. Suppose $x:[t_0,\infty) \to \mathbb{R}$ is a non-negative  continuously differentiable function. Then,  $ \forall t\geq t_0$, $\forall \alpha \in (0, 1]$, the  inequality 
\begin{align}\label{nr4_eq1}
^{C}\!D^{\alpha}_{t_0,t}\left\lbrace{\phi(t)x^{\beta}(t)}\right\rbrace \leq {\phi(t)} \beta{x^{\beta-1}(t)} {^{C}\!D^{\alpha}_{t_0,t}x(t)}, 
\end{align}
holds, if  the inequality 
\begin{align}\label{nr4_eq2}
^{C}\!D^{\alpha}_{t_0,t}\left\lbrace{x^{\beta}(t)}\right\rbrace \leq \beta{x^{\beta-1}(t)} {^{C}\!D^{\alpha}_{t_0,t}x(t)}, 
\end{align}
holds, where the real constant $\beta \geq 1$.
\end{result}
\begin{proof}
We can write
\begin{align}\label{nr4_eq3}
^{C}\!D^{\alpha}_{t_0,t}\left\lbrace{\phi(t)x^{\beta}(t)}\right\rbrace -{\phi(t)} \beta{x^{\beta-1}(t)} {^{C}\!D^{\alpha}_{t_0,t}x(t)} &= {^{C}\!D^{\alpha}_{t_0,t}\left\lbrace{\phi(t)x^{\beta}(t)}\right\rbrace - {\phi(t)} {^{C}\!D^{\alpha}_{t_0,t}x^{\beta}(t)}}  \nonumber \\
&+ \phi(t)\left[ ^{C}\!D^{\alpha}_{t_0,t}\left\lbrace{x^{\beta}(t)}\right\rbrace -\beta{x^{\beta-1}(t)} {^{C}\!D^{\alpha}_{t_0,t}x(t)}\right].
\end{align}
Thus the result follows by using the Result~\ref{NR3} and the inequality \eqref{nr4_eq2} in the equation \eqref{nr4_eq3}.
\end{proof}

\begin{result}
Let $\phi:[t_0,\infty)\to \mathbb{R}$ be  a non-negative, monotonically decreasing and continuously differentiable function. Suppose $x:[t_0,\infty) \to \mathbb{R}$ is a continuously differentiable function. Then,  $ \forall t\geq t_0$, $\forall \alpha \in (0, 1]$, the  inequality 
\begin{align}
^{C}\!D^{\alpha}_{t_0,t}\left\lbrace{\phi(t)x^{\beta}(t)}\right\rbrace \leq {\phi(t)} \beta{x^{\beta-1}(t)} {^{C}\!D^{\alpha}_{t_0,t}x(t)}, 
\end{align}
holds, if  the inequality 
\begin{align}
^{C}\!D^{\alpha}_{t_0,t}\left\lbrace{x^{\beta}(t)}\right\rbrace \leq \beta{x^{\beta-1}(t)} {^{C}\!D^{\alpha}_{t_0,t}x(t)}, 
\end{align}
holds, where the real constant $\beta=2n$, and $n\in \mathbb{N}$.
\end{result}
\begin{proof}
The proof is a consequence of the Result~\ref{NR4}.
\end{proof}

\begin{result}\label{NR5}
Let $\phi:[t_0,\infty)\to \mathbb{R}$ be  a non-negative, monotonically decreasing and continuously differentiable function. Suppose $x:[t_0,\infty) \to \mathbb{R}$ is a positive continuously differentiable function. Then,  $ \forall t\geq t_0$, $\forall \alpha \in (0, 1]$, the  inequality 
\begin{align}\label{nr5_eq1}
^{C}\!D^{\alpha}_{t_0,t}\left\lbrace{\phi(t)x^{\beta}(t)}\right\rbrace \leq {\phi(t)} \beta{x^{\beta-1}(t)} {^{C}\!D^{\alpha}_{t_0,t}x(t)}, 
\end{align}
holds, if  the inequality 
\begin{align}\label{nr5_eq2}
^{C}\!D^{\alpha}_{t_0,t}\left\lbrace{x^{\beta}(t)}\right\rbrace \leq \beta{x^{\beta-1}(t)} {^{C}\!D^{\alpha}_{t_0,t}x(t)}, 
\end{align}
holds, where the real constant $\beta \geq 0$.
\end{result}
\begin{proof}
The proof is similar to the proof of the Result~\ref{NR4}.
\end{proof}
\begin{result}\label{NR6}
Let $\phi:[t_0,\infty)\to \mathbb{R}$ be  a positive, monotonically decreasing and continuously differentiable function. Suppose $x:[t_0,\infty) \to \mathbb{R}$ is a positive continuously differentiable function. Then,  $ \forall t\geq t_0$, $\forall \alpha \in (0, 1]$, the  inequality 
\begin{align}\label{nr6_eq1}
^{C}\!D^{\alpha}_{t_0,t}\left\lbrace{x^{\beta}(t)}\right\rbrace \leq \beta{x^{\beta-1}(t)} {^{C}\!D^{\alpha}_{t_0,t}x(t)}, 
\end{align}
holds, if  the inequality 
\begin{align}\label{nr6_eq2}
^{C}\!D^{\alpha}_{t_0,t}\left\lbrace{\phi(t)x^{\beta}(t)}\right\rbrace \leq {\phi(t)} \beta{x^{\beta-1}(t)} {^{C}\!D^{\alpha}_{t_0,t}x(t)}, 
\end{align}
holds, where the real constant $\beta \geq 0$.
\end{result}
\begin{proof}
Note that 
\begin{align}\label{nr6_eq3}
^{C}\!D^{\alpha}_{t_0,t}\left\lbrace{x^{\beta}(t)}\right\rbrace -\beta{x^{\beta-1}(t)} {^{C}\!D^{\alpha}_{t_0,t}x(t)} &= \frac{1}{\phi(t)}\left[{^{C}\!D^{\alpha}_{t_0,t}\left\lbrace{\phi(t)x^{\beta}(t)}\right\rbrace - {\phi(t)} \beta{x^{\beta-1}(t)}{^{C}\!D^{\alpha}_{t_0,t}x(t)}}\right]  \nonumber \\
&\qquad{}+\frac{1}{\phi(t)}\left[\phi(t){^{C}\!D^{\alpha}_{t_0,t}x^{\beta}(t)}- ^{C}\!D^{\alpha}_{t_0,t}\left\lbrace{\phi(t)x^{\beta}(t)}\right\rbrace \right].
\end{align}
Let us denote by
\begin{align}
f(t)=\frac{1}{\phi(t)}\left[{^{C}\!D^{\alpha}_{t_0,t}\left\lbrace{\phi(t)x^{\beta}(t)}\right\rbrace - {\phi(t)} \beta{x^{\beta-1}(t)}{^{C}\!D^{\alpha}_{t_0,t}x(t)}}\right],
\end{align}
and
\begin{align}
g(t)=\frac{1}{\phi(t)}\left[\phi(t){^{C}\!D^{\alpha}_{t_0,t}x^{\beta}(t)}- {^{C}\!D^{\alpha}_{t_0,t}\left\lbrace{\phi(t)x^{\beta}(t)}\right\rbrace} \right].
\end{align}
Clearly $f(t)\leq 0$, by inequality \eqref{nr6_eq2}.  Let $\psi(t)=\frac{1}{\phi(t)}$. Since $\psi(t)$ is positive, monotonically increasing and continuously differentiable function, it follows from Result~\ref{Nr1}, that $g(t)\leq 0$. This completes the proof.
\end{proof}

\subsection{Main properties}
\begin{lemma}\emph{\cite{Dai2017}}\label{lem3}
Suppose $x:[t_0,\infty) \to \mathbb{R}$ is a non-negative continuously differentiable function. Then,  $ \forall t\geq t_0$, $\forall \alpha \in (0, 1)$, the following inequality holds
\begin{align}
^{C}\!D^{\alpha}_{t_0,t}\left\lbrace{x^{\beta}(t)}\right\rbrace \leq \beta{x^{\beta-1}(t)} {^{C}\!D^{\alpha}_{t_0,t}x(t)}, 
\end{align}
where $\beta \geq 1 $.
\end{lemma}

\begin{lemma}\emph{\cite{Dai2017}}\label{lem4}
Suppose $x:[t_0,\infty) \to \mathbb{R}$ is a continuously differentiable function. Then,  $ \forall t\geq t_0$, $\forall \alpha \in (0, 1)$, the following inequality holds
\begin{align}
^{C}\!D^{\alpha}_{t_0,t}\left\lbrace{x^{\beta}(t)}\right\rbrace \leq \beta{x^{\beta-1}(t)} {^{C}\!D^{\alpha}_{t_0,t}x(t)}, 
\end{align}
where $\beta =\frac{p}{q}\geq 1$, $p=2n$, and $n, q\in \mathbb{N}$.
\end{lemma}

\begin{result}\label{NR7}
Let $\phi:[t_0,\infty)\to \mathbb{R}$ be  a non-negative, monotonically decreasing and continuously differentiable function. Suppose $x:[t_0,\infty) \to \mathbb{R}$ is a continuously differentiable function. Then,  $ \forall t\geq t_0$, $\forall \alpha \in (0, 1]$, the following inequality holds
\begin{align}\label{NR7_eq1}
^{C}\!D^{\alpha}_{t_0,t}\left\lbrace{\phi(t)x^{\beta}(t)}\right\rbrace \leq {\phi(t)} \beta{x^{\beta-1}(t)} {^{C}\!D^{\alpha}_{t_0,t}x(t)}, 
\end{align}
where $\beta =\frac{p}{q}\geq 1$, $p=2k$,  and $k, q\in \mathbb{N}$.
\end{result}
\begin{proof}
The proof follows by using the Lemma~\ref{lem4}. Also it can be easily observed that the inequality \eqref{NR7_eq1} holds from the Result~\ref{NR4} where the Lemma~\ref{lem4} is used.
\end{proof}

\begin{result}\label{NR8}
Let $\phi:[t_0,\infty)\to \mathbb{R}$ be  a non-negative, monotonically decreasing and continuously differentiable function. Suppose $x:[t_0,\infty) \to \mathbb{R}$ is a non-negative continuously differentiable function. Then,  $ \forall t\geq t_0$, $\forall \alpha \in (0, 1]$, the following inequality holds
\begin{align}\label{NR8_eq1}
^{C}\!D^{\alpha}_{t_0,t}\left\lbrace{\phi^{p}(t)x^{\beta}(t)}\right\rbrace \leq {\phi^{p}(t)} \beta{x^{\beta-1}(t)} {^{C}\!D^{\alpha}_{t_0,t}x(t)}, 
\end{align}
where the real positive constants $p\geq 1$ and $\beta \geq 1 $.
\end{result}
\begin{proof}
It follows from the Result~\ref{NR3} that
\begin{align}\label{NR8_eq2}
^{C}\!D^{\alpha}_{t_0,t}\left\lbrace{\phi^{p}(t)x^{\beta}(t)}\right\rbrace \leq {\phi^{p}(t)}  {^{C}\!D^{\alpha}_{t_0,t}x^{\beta}(t)}, 
\end{align}
where  the real positive constants $p\geq 1$ and $\beta \geq 1 $. Then, the application of the Result~\ref{NR4} or the Lemma~\ref{lem3} into the inequality \eqref{NR8_eq2} yields the inequality \eqref{NR8_eq1}.
\end{proof}

\begin{result}\label{NR9}
Let $\phi_{i}:[t_0,\infty)\to \mathbb{R}$ are non-negative, monotonically decreasing and continuously differentiable functions for $i=1,2,\cdots,n$. Suppose $x_{i}:[t_0,\infty) \to \mathbb{R}$ are non-negative continuously differentiable functions for $i=1,2,\cdots,n$ . Then,  $ \forall t\geq t_0$, $\forall \alpha \in (0, 1]$, the following inequality holds
\begin{align}
^{C}\!D^{\alpha}_{t_0,t}\left\lbrace{\sum\limits_{i=1}^{n}\phi_{i}^{p_{i}}(t)x_{i}^{\beta_{i}}(t)}\right\rbrace \leq \sum\limits_{i=1}^{n}{\phi_{i}^{p_{i}}(t)} \beta_{i}{x_{i}^{\beta_{i}-1}(t)} {^{C}\!D^{\alpha}_{t_0,t}x_{i}(t)}, 
\end{align}
where the real positive constants $p_{i}\geq 1$ and $\beta_{i} \geq 1 $ for $i=1,2,\cdots,n$.
\end{result}
\begin{proof}
The proof is similar to the proof of Result~\ref{NR8}.
\end{proof}

\begin{result}\label{NR10}
Let $\phi_{i}:[t_0,\infty)\to \mathbb{R}$ are non-negative, monotonically decreasing and continuously differentiable functions for $i=1,2,\cdots,n$. Suppose $x_{i}:[t_0,\infty) \to \mathbb{R}$ are continuously differentiable functions for $i=1,2,\cdots,n$ . Then,  $ \forall t\geq t_0$, $\forall \alpha \in (0, 1]$, the following inequality holds
\begin{align}
^{C}\!D^{\alpha}_{t_0,t}\left\lbrace{\sum\limits_{i=1}^{n}\phi_{i}^{p_{i}}(t)x_{i}^{\beta_{i}}(t)}\right\rbrace \leq \sum\limits_{i=1}^{n}{\phi_{i}^{p_{i}}(t)} \beta_{i}{x_{i}^{\beta_{i}-1}(t)} {^{C}\!D^{\alpha}_{t_0,t}x_{i}(t)}, 
\end{align}
where the real positive constants $p_{i}\geq 1$ and $\beta_{i} =\frac{u_i}{v_i}\geq 1 $, $u_i=2k_{i}$, $k_{i}, v_{i}\in \mathbb{N}$  for $i=1,2,\cdots,n$.
\end{result}
\begin{proof}
The proof is a consequence of the Result~\ref{NR9} where the Lemma~\ref{lem4} is used.
\end{proof}

\begin{result}\label{NR11}
Let $\phi_{i}:[t_0,\infty)\to \mathbb{R}$ are non-negative, monotonically decreasing and continuously differentiable functions for $i=1,2,\cdots,n$. Suppose $x_{i}:[t_0,\infty) \to \mathbb{R}$ are continuously differentiable functions for $i=1,2,\cdots,n$ . Then,  $ \forall t\geq t_0$, $\forall \alpha \in (0, 1]$, the following inequality holds
\begin{align}
^{C}\!D^{\alpha}_{t_0,t}\left\lbrace{\sum\limits_{i=1}^{n}c_{i}\phi_{i}^{p_{i}}(t)x_{i}^{\beta_{i}}(t)+\sum\limits_{i=1}^{n}d_{i}x_{i}^{\beta_{i}}(t)}\right\rbrace &\leq \sum\limits_{i=1}^{n}c_{i}{\phi_{i}^{p_{i}}(t)} \beta_{i}{x_{i}^{\beta_{i}-1}(t)} {^{C}\!D^{\alpha}_{t_0,t}x_{i}(t)}\\
&\qquad{}+\sum\limits_{i=1}^{n}d_{i} \beta_{i}{x_{i}^{\beta_{i}-1}(t)} {^{C}\!D^{\alpha}_{t_0,t}x_{i}(t)}, 
\end{align}
where the real positive constants $c_i>0$,  $d_i>0$, $p_{i}\geq 1$, and $\beta_{i} =\frac{u_i}{v_i}\geq 1 $, $u_i=2k_{i}$, $k_{i}, v_{i}\in \mathbb{N}$  for $i=1,2,\cdots,n$.
\end{result}
\begin{proof}
The result follows by combining the Lemma~\ref{lem4} with the Result~\ref{NR10}.
\end{proof}

\begin{result}\label{NR12}
Let $\phi_{i}:[t_0,\infty)\to \mathbb{R}$ are non-negative, monotonically decreasing and continuously differentiable functions for $i=1,2,\cdots,n$. Suppose $x_{i}:[t_0,\infty) \to \mathbb{R}$ are continuously differentiable functions for $i=1,2,\cdots,n$. Then,  $ \forall t\geq t_0$, $\forall \alpha \in (0, 1]$, the following inequality holds
\begin{align}
^{C}\!D^{\alpha}_{t_0,t}\left\lbrace{\sum\limits_{i=1}^{n}c_{i}\phi_{i}^{p_{i}}(t)x_{i}^{\beta_{i}}(t) +\sum\limits_{i=1}^{n}d_{i}x_{i}^{\beta_{i}}(t)}\right\rbrace &{+} ^{C}\!D^{\alpha}_{t_0,t}\left\lbrace{\sum\limits_{i=1}^{n}a_{i}x_{i}^{\gamma_{i}}(t)}\right\rbrace   \nonumber\\
 & \leq\sum\limits_{i=1}^{n}c_{i}{\phi_{i}^{p_{i}}(t)} \beta_{i}{x_{i}^{\beta_{i}-1}(t)} {^{C}\!D^{\alpha}_{t_0,t}x_{i}(t)} \nonumber\\
 &\qquad{}+\sum\limits_{i=1}^{n}d_{i} \beta_{i}{x_{i}^{\beta_{i}-1}(t)} {^{C}\!D^{\alpha}_{t_0,t}x_{i}(t)} \nonumber\\
 &\qquad{}+\sum\limits_{i=1}^{n}a_{i} \gamma_{i}{x_{i}^{\gamma_{i}-1}(t)} {^{C}\!D^{\alpha}_{t_0,t}x_{i}(t)},
\end{align}
where the real constants $a_i\geq 0$, $c_i\geq0$,  $d_i\geq0$, $p_{i}\geq 1$,  $\beta_{i} =\frac{u_i}{v_i}\geq 1 $, $u_i=2k_{i}$, $k_{i}, v_{i}\in \mathbb{N}$,  and $\gamma_{i} =\frac{m_i}{n_i}\geq 1 $,  $m_i=2\ell_{i}$, $\ell_{i}, n_{i}\in \mathbb{N}$  for $i=1,2,\cdots,n$.
\end{result}
\begin{proof}
The proof is a consequence of the Result~\ref{NR11}.
\end{proof}

\section{Examples}
In this section, we discuss the fractional order generalizations or versions of a few interesting examples \cite{Khalil2002}. We apply the fractional Lyapunov direct method and utilize appropriate fractional differential inequalities, in order to analyse the stability of the zero solutions to nonautonomous fractional order systems. We use the numerical predictor-corrector method \cite{Diethelm2002} in order to carry out the numerical solutions to the presented examples.

\begin{example}\label{MEx2}
Consider the  nonautonomous linear fractional order system
\begin{equation}\label{MEx2_eq1}
\begin{aligned}
^{C}\!D^{\alpha}_{0,t}x_{1}(t)&=-x_{1}(t)-\frac{1}{1+t}x_{2}(t),\;\; x_{1}(0)=x_{10}, \\
^{C}\!D^{\alpha}_{0,t}x_{2}(t)&=x_{1}(t)-x_{2}(t),\;\; x_{2}(0)=x_{20},
\end{aligned}
\end{equation}
where $0<\alpha\leq 1$.
\end{example}

Let us choose the function $V(t,x)=x_{1}^{2}+x_{2}^{2}+\frac{1}{1+t}x_{2}^{2}$. Note that $x_{1}^{2}+x_{2}^{2}\leq V(t,x)\leq x_{1}^{2}+2x_{2}^{2}$, $\forall {x=\left(x_1,x_2\right)^{T}\in \mathbb{R}^{2}}$. Then, the application of the Result~\ref{NR12}, allows us to calculate the Caputo fractional derivative of $V(t,x)$ along the solution $x(t)$ to \eqref{MEx2_eq1} as follows 
\begin{align}\label{MEx2_eq2}
^{C}\!D^{\alpha}_{0,t}V(t,x(t)) &\leq \left[-2x_{1}^{2}(t)-\frac{2}{1+t}x_{1}(t)x_{2}(t)\right]+\left[(2+\frac{2}{1+t})x_{1}(t)x_{2}(t)-(2+\frac{2}{1+t})x_{2}^{2}(t)\right]  \nonumber \\
& =-2x_{1}^{2}(t)+2x_{1}(t)x_{2}(t)-(2+\frac{2}{1+t})x_{2}^{2}(t) \nonumber  \\
& \leq -2x_{1}^{2}(t)+2x_{1}(t)x_{2}(t)-2x_{2}^{2}(t)  \nonumber  \\
&=-x^{T}(t)\left(\begin{array}{cc}
2 & -1\\
-1 & 2
\end{array}\right) x(t), \;\; \forall\; {x=\begin{pmatrix}
x_1\\x_2
\end{pmatrix}\in \mathbb{R}^{2}}.
\end{align}
Note that $x^{T}\left(\begin{array}{cc}
2 & -1\\
-1 & 2
\end{array}\right)x$ is a positive definite quadratic function. Then, it follows from \eqref{MEx2_eq2} that 
 \begin{equation}\label{MEx2_eq3}
^{C}\!D^{\alpha}_{0,t}V(t,x(t)) \leq -\Vert{x(t)}\Vert^{2}.
\end{equation}
Let $\gamma_1(r)=r^{2}$, $\gamma_2(r)=2r^{2}$ and $\gamma_3(r)=r^{2}$, where $r=\Vert{x}\Vert$. Then, all the assumptions of Theorem 6.2 of \cite{YLi2010} are satisfied. This observation confirms that $V(t,x)$ is indeed a time-varying Lyapunov function.  Hence, by Theorem 6.2 of  \cite{YLi2010}, we conclude that the zero solution is asymptotically stable.

\begin{figure}[H]
\centering
\includegraphics[width=4.5in,height=3in]{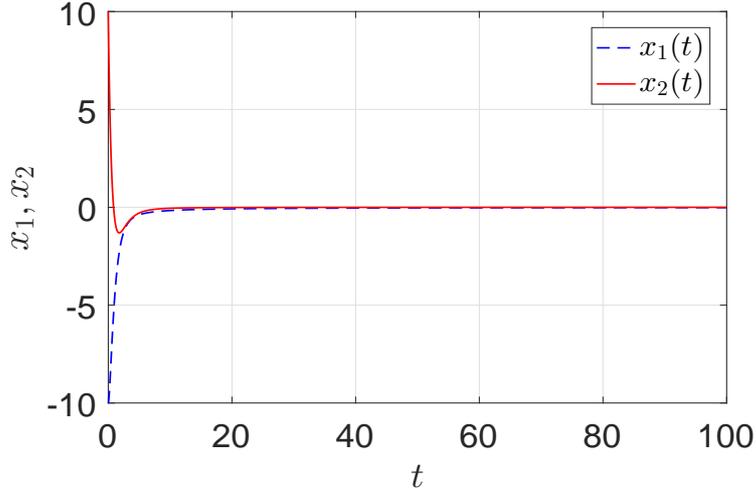}
\caption{Numerical solution to \eqref{MEx2_eq1} for the value of fractional order $\alpha=0.9$ starting from initial values $x_{1}(0)=-10$ and $x_{2}(0)=10$.}\label{MEx2_fig1}
\end{figure}

Further, from \eqref{MEx2_eq3}, we deduce
\begin{equation}
V(t,x(t))\leq E_{\alpha}\left(-\frac{1}{2}t^{\alpha}\right)V(0,x(0)),\;\; \forall {t\geq 0}.
\end{equation}
Thus, it follows that 
\begin{equation}
\Vert{x(t)}\Vert\leq \sqrt{2E_{\alpha}\left(-\frac{1}{2}t^{\alpha}\right)\Vert{x(0)}\Vert ^{2}},\;\; \forall {t\geq 0}.
\end{equation}
As a result, the zero solution to the system \eqref{MEx2_eq1} is Mittag-Leffler stable. Thus, the zero solution is asymptotically stable. The numerical solution shown in the Figure~\ref{MEx2_fig1} indicates that the non-trivial solutions approach to the zero solution. 

\begin{figure}[H]
\centering
\includegraphics[width=4.5in,height=3in]{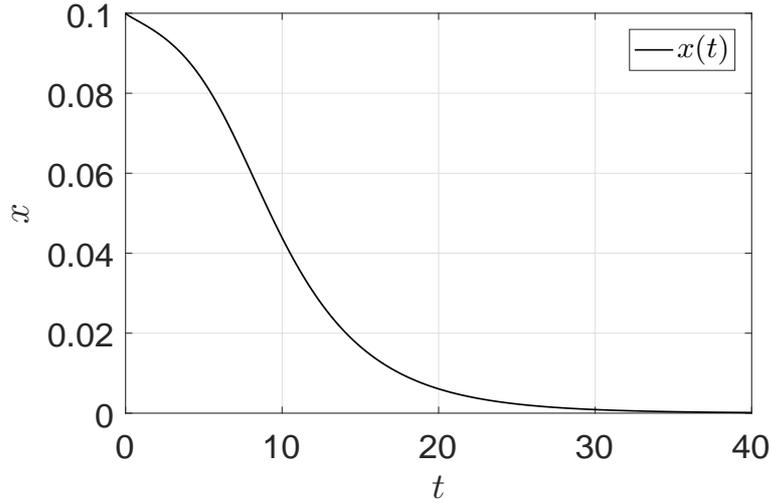}
\caption{Numerical solution to \eqref{Ex1_eq1} for the value of fractional order $\alpha=0.8$ with initial value $x(0)=0.1$.}\label{Ex1_fig1}
\end{figure}

\begin{example}\label{Ex1}
Consider the nonautonomous nonlinear fractional differential equation
\begin{align}\label{Ex1_eq1}
^{C}\!D^{\alpha}_{0,t}x(t)=-x^{3}(t)-e^{t/2}x^{3}(t),\;\; x(0)=x_{0}, \;\; 0<\alpha\leq 1.
\end{align}
\end{example}
Let $V(t,x)=x^{6}+ e^{-t/2}x^{6}$ be the function, which depends on time $t$ and variable $x$. Then, by using the  Result~\ref{NR12}, we get the Caputo fractional derivative of $V(t,x)$ along the solution $x(t)$ to \eqref{Ex1} as follows
\begin{align}
^{C}\!D^{\alpha}_{0,t}V(t,x(t)) &\leq \left[-6x^{8}(t)-6e^{t/2}x^{8}(t)\right] +\left[-6e^{-t/2}x^{8}(t)-6x^{8}(t)\right] \nonumber \\
& =-6\left(1+e^{-t/2}\right)x^{8}(t) -6\left(1+e^{t/2}\right)x^{8}(t)  \nonumber \\
& \leq -12x^{8}(t), \;\; \forall\; {x\in \mathbb{R}},\;\; \forall\; {t\geq{0}}.
\end{align}
Note that $x^{6}\leq V(t,x)\leq 2x^{6}$, $\forall {x\in \mathbb{R}}$, $\forall {t\geq 0}$.  Let $\gamma_1(r)=r^{6}$, $\gamma_2(r)=2r^{6}$ and $\gamma_3(r)=12r^{8}$, where $r=\vert{x}\vert$. Then, we see that all the assumptions of  the Theorem 6.2 of  \cite{YLi2010} are satisfied. Hence, it follows from the Theorem 6.2 of \cite{YLi2010} that the zero solution to the equation \eqref{Ex1_eq1} is asymptotically stable. 
The numerical solution is shown in the Figure~\ref{Ex1_fig1}.

\begin{example}\label{Ex3}
Consider the nonautonomous nonlinear fractional order system 
\begin{equation}\label{Ex3_eq1}
\begin{aligned}
^{C}\!D^{\alpha}_{0,t}x_{1}(t)&=-x_{1}(t)-x_{2}(t)+\sin(t)\left(x_{1}^{2}(t)+x_{2}^{2}(t)\right),\;\; x_{1}(0)=x_{10}, \\
^{C}\!D^{\alpha}_{0,t}x_{2}(t)&=x_{1}(t)-x_{2}(t)+\cos(t)\left(x_{1}^{2}(t)+x_{2}^{2}(t)\right),\;\; x_{2}(0)=x_{20},
\end{aligned}
\end{equation}
where $0<\alpha\leq 1$.

\end{example}

Let $V(t,x)=\frac{x_{1}^{2}}{2}+ \frac{x_{2}^{2}}{2}+\phi(t)\frac{x_{1}^{2}}{2}+\phi(t)\frac{x_{2}^{2}}{2}$ be the time-varying Lyapunov function, where $\phi(t)$ is a non-negative, monotonically decreasing, bounded and continuously differentiable function. Then, by using the  Result~\ref{NR12} (or the Result~\ref{Nr3} where the Lemma 1 of \cite{Camacho2014} is used), we get the Caputo fractional derivative of $V(t,x)$ along the solution $x(t)$ to \eqref{Ex3} as follows
\begin{align}
^{C}\!D^{\alpha}_{0,t}V(t,x(t)) &\leq \left[-x_{1}^{2}(t)-x_{1}(t)x_{2}(t)+x_{1}(t)\sin(t)\left(x_{1}^{2}(t)+x_{2}^{2}(t)\right)\right] \nonumber\\
&\quad{}+\left[x_{1}(t)x_{2}(t)-x_{2}^{2}(t)+x_{2}(t)\cos(t)\left(x_{1}^{2}(t)+x_{2}^{2}(t)\right)\right] \nonumber \\
&\quad{}+\phi(t)\left[-x_{1}^{2}(t)-x_{1}(t)x_{2}(t)+x_{1}(t)\sin(t)\left(x_{1}^{2}(t)+x_{2}^{2}(t)\right)\right] \nonumber\\
&\quad{}+\phi(t)\left[x_{1}(t)x_{2}(t)-x_{2}^{2}(t)+x_{2}(t)\cos(t)\left(x_{1}^{2}(t)+x_{2}^{2}(t)\right)\right] \nonumber\\
& =-\left[(1+\phi(t))(x_{1}^{2}(t)+x_{2}^{2}(t))\right] \nonumber\\
&\quad{}+\left[(1+\phi(t))(x_{1}^{2}(t)+x_{2}^{2}(t))(x_{1}(t)\sin(t)+x_{2}(t)cos(t))\right]  \nonumber\\
& \leq -\left[(1+\phi(t))(x_{1}^{2}(t)+x_{2}^{2}(t))\right] +\left[(1+\phi(t))(x_{1}^{2}(t)+x_{2}^{2}(t))^{3}\right]  \nonumber\\
& = -(1+\phi(t))\Vert{x(t)}\Vert^{2} (1-\Vert{x(t)}\Vert) \nonumber\\
& \leq -(1+\phi(t))(1-r)\Vert{x(t)}\Vert^{2},\;\; \forall  \Vert{x(t)}\Vert \leq r, \;\; \forall r<1,  \nonumber\\
& \leq -(1-r)\Vert{x(t)}\Vert^{2},\;\; \forall  \Vert{x(t)}\Vert \leq r, \;\; \forall r<1. \nonumber\\
\end{align}
Hence, it follows from the Theorem 6.2 of \cite{YLi2010} that the zero solution to the system \eqref{Ex3_eq1} is asymptotically stable. 
\begin{figure}[H]
\centering
\includegraphics[width=4.5in,height=3in]{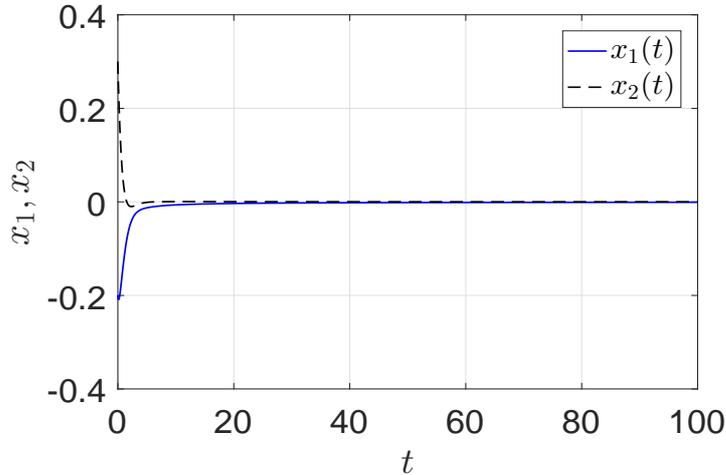}
\caption{Numerical solution to \eqref{Ex3_eq1}, starting from the initial values $x_{1}(0)=-0.2$ and $x_{2}(0)=0.3$, for the value of fractional order $\alpha=0.85$.}\label{Ex3_fig1}
\end{figure}
The numerical solution to the system \eqref{Ex3_eq1} is shown in the Figure~\ref{Ex3_fig1} for some suitable value of fractional order $\alpha=0.85$. It indicates that the zero solution is asymptotically stable.

\section{Conclusions}
We have established several new fractional differential inequalities which involve the Caputo fractional derivative of continuously differentiable functions. By presenting a few simple illustrative examples, applying fractional direct method and using appropriate fractional differential inequalities, we have shown that the presented fractional differential inequalities play a basic role in the  analysis as well as the calculation of bounds of solutions to nonautonomous fractional order systems. These fractional differential inequalities open up an opportunity and show a way to construct or discover different types of potential candidate time-varying Lyapunov functions. They can also be useful for the calculation of their fractional derivatives, at the same time, for the stability analysis of fractional order systems.


\begin{thebibliography}{30}

\bibitem{Podlubny1999}
 I. Podlubny, Fractional differential equations, San Diego,  Academic Press, 1999.

\bibitem{Kilbas2006}
 A. A. Kilbas, H. M. Srivastava, J. J. Trujillo, Theory and applications of fractional differential equations, Amsterdam,  Elsevier, 2006.

\bibitem{Petras2011}
I. Petr\'{a}\v{s}, Fractional order nonlinear systems: Modelling, analysis and simulation, Springer, 2010.


\bibitem{Kiryakova2014}
D. Val\'{e}rio, J. T. Machado, V. Kiryakova, Some pioneers of the applications of fractional calculus, Fractional Calculus and Applied Analysis, 17, (2014), pp. 552--578.

\bibitem{Tarasov2013}
V. E. Tarasov, No violation of the Leibniz rule. No fractional derivative, Communications in Nonlinear Science and Numerical Simulation, 18, (2013), pp. 2945--2948.

\bibitem{VETarasov2016}
V. E. Tarasov, Leibniz rule and fractional derivatives of power functions, Journal of Computational and Nonlinear Dynamics, 11, (2016), 031014--4.

\bibitem{Tarasov2016}
V. E. Tarasov, On chain rule for fractional derivatives,  Communications in Nonlinear Science and Numerical Simulation,   30, (2016), pp. 1--4.


\bibitem{YLi2009}
Y. Li, Y. Q. Chen, I. Podlubny, Mittag-leffler stability of fractional order nonlinear dynamic systems, Automatica, 45, (2009), pp. 1965--1969.

\bibitem{YLi2010}
Y. Li, Y. Q. Chen, I. Podlubny, Stability of fractional-order nonlinear dynamic systems: Lyapunov direct method and generalized Mittag-Leffler stability, Computer Mathematics and Applications,  59, (2010), pp. 1810--1821.

\bibitem{Delavari2012}
H. Delavari, D. Baleanu, J. Sadati, Stability analysis of caputo fractional-order nonlinear systems revisited, Nonlinear Dynamics,  67, (2012), pp. 2433--2439 .

\bibitem{Yu2013}
J. Yu, H. Hu, S. Zhou, X. Lin, Generalized Mittag-Leffler stability of multi-variables fractional order nonlinear systems,  Automatica, 49, (2013), pp. 1798--803 .


\bibitem{Mermoud2015}
M. A. Duarte-Mermoud, N. Aguila-Camacho, J. A. Gallegos, R. Castro-Linares, Using general quadratic Lyapunov functions to prove Lyapunov uniform stability for fractional order systems, Communications in Nonlinear Science and Numerical Simulation,   22, (2015), pp. 650--659.
 

\bibitem{BKL2019}
B. K. Lenka, Fractional comparison method and asymptotic stability of multivariable fractional order systems, Communications in Nonlinear Science and Numerical Simulation, 69, (2019), pp. 398--415.


\bibitem{Fernandez2017}
G. Fern\'{a}ndez-Anayaa, G. Nava-Antonio, J. Jamous-Galante, R. Mu\~{n}oz-Vega, E. G. Hern\'{a}ndez-Mart\'{i}nez, Asymptotic stability of distributed order nonlinear dynamical systems,  Communications in Nonlinear Science and Numerical Simulation,   48, (2017), pp. 541--549.


\bibitem{Khalil2002}
H. K. Khalil, Nonlinear systems, Prentice Hall, 2002.

\bibitem{Camacho2014}
N. Aguila-Camacho, M. A. Duarte-Mermoud, J. A. Gallegos, Lyapunov functions for fractional order systems, Communications in Nonlinear Science and Numerical Simulation, 19, (2014), pp. 2951--2957.


\bibitem{Dai2017}
H. Dai, W. Chen, New power law inequalities for fractional derivative and stability analysis of fractional order systems,  Nonlinear Dynamics, 87, (2017), pp. 1531--1542.

\bibitem{Fernandez2017b}
G. Fern\'{a}ndez-Anayaa, G. Nava-Antonio, J. Jamous-Galante, R. Mu\~{n}oz-Vega, E. G. Hern\'{a}ndez-Mart\'{i}nez, Lyapunov functions for a class of nonlinear systems using Caputo derivative,  Communications in Nonlinear Science and Numerical Simulation, 43, (2017), pp. 91--99.

\bibitem{Fernandez2018}
G. Fern\'{a}ndez-Anayaa, G. Nava-Antonio, J. Jamous-Galante, R. Mu\~{n}oz-Vega, E. G. Hern\'{a}ndez-Mart\'{i}nez, Corrigendum to ``Lyapunov functions for a class of nonlinear systems using Caputo derivative" [Commun Nonlinear Sci Numer Simulat 43 (2017) 91--99]. Communications in Nonlinear Science and Numerical Simulation, 56, (2018), pp. 596--597.

\bibitem{Kirane2015}
A. Alsaedi, B. Ahmad, M. Kirane, Maximum principle for certain generalized time and space fractional diffusion equations, Quarterly of Applied Mathematics, 73 (2015), pp. 163--175. 

\bibitem{Kirane2017}
A. Alsaedi, B. Ahmad, M. Kirane, A survey of useful inequalities in fractional calculus, Fractional Calculus and Applied Analysis, 20, (2017), pp. 574--594.

\bibitem{BKLenka2019}
B. K. Lenka, Time-varying Lyapunov functions and  Lyapunov
stability of nonautonomous fractional order systems, International Journal of Applied Mathematics, 32, (2019), pp. 111--130.

\bibitem{Diethelm2002}
K. Diethelm, N. J. Ford, A. D. Freed. A predictor-corrector approach for the numerical solution of fractional differential equations. Nonlinear Dynamics, 2002, 29, pp. 3--22.
\end{thebibliography}
\end{document}